\documentclass[number]{amsart}
\usepackage{amsmath, latexsym, amssymb, amsfonts,amscd,enumerate}
\usepackage{epstopdf}
\usepackage{graphicx}
\usepackage[numbers]{natbib}
\usepackage[myheadings]{fullpage}
\usepackage{epstopdf}
\usepackage{subfigure}
\numberwithin{equation}{section}

\title{Combinatorial Proofs of an Identity from Ramanujan's Lost Notebook and its Variations}

\author{Paul Levande \\ plevande@math.upenn.edu}
\address{Department of Mathematics, University of Pennsylvania, David Rittenhouse Lab, 209 South 33rd Street, Philadelphia, PA 19104-6395}
\begin{document}
\begin{abstract}We examine an identity originally stated in Ramanujan's ``lost notebook'' and first proven algebraically by Andrews and combinatorially by Kim.  We give two independent combinatorial proofs and interpretations of this identity, which also extends an identity recently proven by Pak and Waarnar related to the product of partial theta functions: First, we give a direct combinatorial proof, using the involution principle, of a special case of the identity, and extend this into a direct combinatorial proof of the full identity as written.  Second, we show that the identity can be rewritten, using minor algebraic manipulation, into an identity that can be proven with a direct bijection.  We provide such a bijection using a generalization of a standard bijection from partition theory.     
\end{abstract}
 \maketitle

\pagestyle{myheadings}
\markboth{Paul Levande}{Combinatorial Ramanujan Proofs}
\section{Introduction} In his 1979 introduction to Ramanujan's ``lost'' notebook, Andrews \cite{Andrewsint} gave the following identity, first found there:
\begin{eqnarray} \label{ramorigin}
1+ \sum_{n=1}^{\infty} \frac{q^{n}}{(1-aq)(1-aq^{2}) \cdots (1-aq^{n})(1-bq)(1-bq^{2}) \cdots (1-bq^{n})}   \\
 = (1-a^{-1}) \left(1+ \sum_{n=1}^{\infty} \frac{(-1)^{n} q^{\binom{n+1}{2}} b^{n}a^{-n}} {(1-bq)(1-bq^{2}) \cdots (1-bq^{n})} \right)  \nonumber \\ 
+ a^{-1}\sum_{n=0}^{\infty} (-1)^{n} q^{\binom{n+1}{2}}b^{n}a^{-n} \prod_{k=1}^{\infty} \frac{1}{(1-aq^{k})(1-bq^{k})}  \nonumber
\end{eqnarray}
as well as a proof using algebraic manipulation and an identity of Rogers' rather than a combinatorial interpretation of the identity.  In his recent survey of partition bijections, Pak listed \cite{paksurvey} finding a combinatorial proof of (\ref{ramorigin}) as an unsolved problem.  (Our statement of the identity follows Pak's notation).  Pak has also noted \cite{pakinv} as a general principle that partition identities apparently requiring a combinatorial proof based on a sign-reversing involution, like (\ref{ramorigin}), can often be interpreted as special cases of partition identities that can be proven using direct bijections.  

We will show how Ramanujan's identity (\ref{ramorigin}) fits, broadly speaking, into this framework: First, we will provide a direct combinatorial proof, based on the involution principle, of the special case of (\ref{ramorigin}) found by evaluating (\ref{ramorigin}) at $a=1$.  We will expand this into a direct combinatorial proof of (\ref{ramorigin}) in full.  Next, using some minor algebraic manipulation, we will rewrite (\ref{ramorigin}) into a form suggesting a possible direct bijective proof: 

\begin{eqnarray} \label{rambijone}   \sum_{n=0}^{\infty} \frac{aq^{n}\prod_{k=n+1}^{\infty}(1+acq^{k})}{(1-a)(1-aq)(1-aq^{2}) \cdots (1-aq^{n})} + \sum_{n=1}^{\infty} q^{\binom{n+1}{2}} c^{n}\prod_{k=n+1}^{\infty}(1+acq^{k})  
 = \left(\sum_{n=0}^{\infty}q^{\binom{n+1}{2}}c^{n}\right) \left(\prod_{k=0}^{\infty} \frac{1}{(1-aq^{k})} \right)  \end{eqnarray}
We will prove (\ref{rambijone}) using an explicit bijection.  In fact, this bijection will be seen to be, itself, a generalization of a standard bijection from partition theory.

Note: Kim recently and independently independently found \cite{kimpartialtheta} an alternative combinatorial interpretation and proof of this identity.  The same article notes a recent identity also proven \cite{Andrewspartialtheta} recently and independently by Andrews and Warnaar,
\begin{eqnarray*}
\sum_{n=0}^{\infty} (-a)^{n}q^{\frac{(n)(n-1)}{2}} = (a)_{\infty}(q)_{\infty}  \sum_{n=0}^{\infty} \frac{q^{n}}{(a)_{n}(q)_{n}}
\end{eqnarray*}
using the standard notation of $(p)_{n} = (1-p)(1-pq) \ldots (1-pq^{n-1})$ for any $p$.  This identity can be seen as a special case of (\ref{rambijone}), derived by multiplying both sides of (\ref{rambijone}) by $1-a$ and evaluating at $a=1$, multiplying both sides by $(q)_{\infty} = \prod_{i=1}^{\infty} (1-q^{i})$, and replacing $c$ by $-a$.  
\section{A Combinatorial Proof of The Special Case of Ramanujan's Identity at $a=1$}
First, for a partition $A$, let $\ell(A)$ be the number of parts of $A$, $|A|$ be the sum of the parts of $A$, $a(A)$ to be the largest part of $A$, and $s(A)$ be the smallest part of $A$, as usual.  
Evaluating (\ref{ramorigin}) at $a=1$ gives the following identity, which we will now prove using the involution principle:
\begin{eqnarray} \label{ramoriginspecial}  1+ \sum_{n=1}^{\infty} \frac{q^{n}}{(1-q)(1-q^{2}) \cdots (1-q^{n})(1-bq)(1-bq^{2}) \cdots (1-bq^{n})}  \\
 = \sum_{n=0}^{\infty} (-1)^{n} q^{\binom{n+1}{2}}b^{n} \prod_{k=1}^{\infty} \frac{1}{(1-q^{k})(1-bq^{k})} \nonumber 
 \end{eqnarray}  
 Let $S$ be the set of ordered triples $(\lambda, \mu, \gamma)$ such that $\lambda$ is a (possibly empty) ``triangular'' partition--i.e., $\lambda = (n, n-1, \ldots, 1)$ for some $n$--and $\mu$ and $\gamma$ are partitions with positive parts.  Let $S^{+}$ be the subset of $S$ consisting of ordered triples $(\lambda, \mu, \gamma)$ such that $\ell(\lambda)$ is even, with $S^{-}$ the subset of $S$ consisting of ordered triples $(\lambda, \mu, \gamma)$ such that $\ell(\lambda)$ is odd.  
 
Define $\psi: S \rightarrow S$ as follows:
\begin{enumerate}
\item If $\lambda = \emptyset$ and $a(\mu) \geq a(\gamma)$, let $\psi((\lambda, \mu, \gamma)) = (\lambda, \mu, \gamma)$.  
\item If $a(\lambda) + a(\mu) \geq a(\gamma)$ and $a(\lambda) \neq 0$ (i.e., if $\lambda \neq \emptyset$) let $\psi((\lambda, \mu, \gamma)) = (\tilde{\lambda}, \tilde{\mu}, \tilde{\gamma})$, where $\tilde{\lambda} = (\lambda_{2}, \lambda_{3}, \ldots)$, $\tilde{\mu} = (\mu_{2}, \mu_{3}, \ldots)$ and $\tilde{\gamma} = (\lambda_{1}+\mu_{1}, \gamma_{1}, \gamma_{2}, \ldots)$.  
\item If $a(\lambda) + a(\mu) < a(\gamma)$, let $\psi((\lambda, \mu, \gamma)) = (\tilde{\lambda}, \tilde{\mu}, \tilde{\gamma})$, where $\tilde{\lambda} = (\lambda_{1}+1, \lambda_{1}, \lambda_{2}, \lambda_{3}, \ldots)$, $\tilde{\mu} = (\gamma_{1}-\lambda_{1}-1, \mu_{1}, \mu_{2}, \ldots)$ and $\tilde{\gamma} = (\gamma_{2}, \gamma_{3} \ldots)$.  Note that, because $\gamma_{1} - \lambda_{1} > \mu_{1}$, $\gamma_{1} - \lambda_{1} - 1 \geq \mu_{1}$.  
\end{enumerate}

Note that $\lambda_{1} + \mu_{1}> \lambda_{2} + \mu_{2}$, since $\lambda_{1} = \lambda_{2}+1$, and that, therefore, $\psi$ is involutive.  Note further that, in all cases, if $\psi((\lambda, \mu, \gamma)) = (\tilde{\lambda}, \tilde{\mu}, \tilde{\gamma})$, $|\lambda|+|\mu|+|\gamma| = |\tilde{\lambda}|+|\tilde{\mu}|+|\tilde{\gamma}|$ and $\ell(\lambda)+\ell(\gamma) = \ell(\tilde{\lambda}) + \ell(\tilde{\gamma})$.  In addition, either $(\lambda, \mu, \gamma) \in Fix(\psi)$ or $\ell(\lambda) = \ell(\tilde{\lambda}) \pm 1$.  Finally, $(\lambda, \mu, \gamma) \in Fix(\psi)$ if and only if $\lambda = \emptyset$ and $a(\mu) \geq a(\gamma)$, and therefore

\begin{eqnarray*} \sum_{(\lambda, \mu, \gamma) \in Fix(\psi)} q^{|\lambda| + |\mu| + |\gamma|}b^{\ell(\lambda) + \ell(\gamma)} = 1+ \sum_{n=1}^{\infty} \frac{q^{n}}{(1-q) \cdots (1-q^{n})(1-bq) \cdots (1-bq^{n})} \end{eqnarray*}
Since 

\begin{eqnarray*} \sum_{(\lambda, \mu, \gamma) \in S} (-1)^{\ell(\lambda)} q^{|\lambda| + |\mu| + |\gamma|}b^{\ell(\lambda) + \ell(\gamma)} = \sum_{n=0}^{\infty} (-1)^{n} q^{\binom{n+1}{2}}b^{n} \prod_{k=1}^{\infty} \frac{1}{(1-q^{k})(1-bq^{k})} \end{eqnarray*}
this suffices to prove (\ref{ramoriginspecial}) using the involution principle.

\section{A Combinatorial Proof of Ramanujan's Identity As Written}
Let us further examine the above proof of (\ref{ramoriginspecial}) using the involution principle: The involution $\psi$ splits $S$ into $Fix(\psi)$ and unordered pairs $\left\{ (\lambda, \mu, \gamma), (\tilde{\lambda}, \tilde{\mu}, \tilde{\gamma}) \right\}$ of distinct elements of $S$, with $\psi((\lambda, \mu, \gamma)) = (\tilde{\lambda}, \tilde{\mu}, \tilde{\gamma})$, $|\lambda|+|\mu|+|\gamma| = |\tilde{\lambda}|+|\tilde{\mu}|+|\tilde{\gamma}|$, $\ell(\lambda) + \ell(\mu) = \ell(\tilde{\lambda}) + \ell(\tilde{\mu})$, and $\ell(\tilde{\lambda}) = \ell(\lambda) \pm 1$.  Assume we can define disjoint subsets $B_{1}, B_{2}$ of $S$ such that $\psi(B_{1}) = B_{2}$ and some statistic $f$ on $S$ such that $f((\lambda, \mu, \gamma)) = f(\psi((\lambda, \mu, \gamma)))+1$ for all $(\lambda, \mu, \gamma) \in B_{1}$.  Then
\begin{eqnarray*} && \sum_{(\lambda, \mu, \gamma) \in S} (-1)^{\ell(\lambda)} a^{f((\lambda, \mu, \gamma))}b^{\ell(\lambda)+\ell(\gamma)}q^{|\lambda|+|\mu|+|\gamma|} \\ &&= (1-a) \left( \sum_{(\lambda, \mu, \gamma) \in B_{1}}(-1)^{\ell(\lambda)} a^{f((\lambda, \mu, \gamma))}b^{\ell(\lambda)+\ell(\gamma)}q^{|\lambda|+|\mu|+|\gamma|}\right) +  \sum_{(\lambda, \mu, \gamma) \in Fix(\psi)} a^{f((\lambda, \mu, \gamma))}b^{\ell(\lambda)+\ell(\gamma)}q^{|\lambda|+|\mu|+|\gamma|}  \end{eqnarray*}
or 
\begin{eqnarray*} && \sum_{(\lambda, \mu, \gamma) \in S} (-1)^{\ell(\lambda)} a^{f((\lambda, \mu, \gamma))}b^{\ell(\lambda)+\ell(\gamma)}q^{|\lambda|+|\mu|+|\gamma|} - (1-a) \left( \sum_{(\lambda, \mu, \gamma) \in B_{1}}(-1)^{\ell(\lambda)} a^{f((\lambda, \mu, \gamma))}b^{\ell(\lambda)+\ell(\gamma)}q^{|\lambda|+|\mu|+|\gamma|}\right) \\ &&= \sum_{(\lambda, \mu, \gamma) \in Fix(\psi)} a^{f((\lambda, \mu, \gamma))}b^{\ell(\lambda)+\ell(\gamma)}q^{|\lambda|+|\mu|+|\gamma|}  \end{eqnarray*}.
Informally, one can think of the term 
\begin{eqnarray*} (1-a) \left( \sum_{(\lambda, \mu, \gamma) \in B_{1}}(-1)^{\ell(\lambda)} a^{f((\lambda, \mu, \gamma))}b^{\ell(\lambda)+\ell(\gamma)}q^{|\lambda|+|\mu|+|\gamma|}\right) \end{eqnarray*} 
as an ``error term'' cancelling out pairs $\left\{ (\lambda, \mu, \gamma), (\tilde{\lambda}, \tilde{\mu}, \tilde{\gamma}) \right\}$ with $(\lambda, \mu, \gamma) \in B_{1}$ and $(\tilde{\lambda}, \tilde{\mu}, \tilde{\gamma}) \in B_{2}$--pairs that do not quite cancel out in the main alternating sum, since, although the exponents of $q$ and $b$ match, the exponents of $a$ differ by $1$.

Define disjoint subsets $B_{1}, B_{2}$ of $S$ as follows:
\begin{displaymath}  B_{1} = \left\{(\lambda, \mu, \gamma) \in S: \mu = \emptyset, 0 < a(\gamma) \leq a(\lambda) \right\}; B_{2} = \left\{ (\lambda, \mu, \gamma) \in S : \mu = \emptyset; a(\gamma) = a(\lambda)+1 \right\} \end{displaymath}.

Note that $\psi(B_{1}) = B_{2}$, and that if $(\lambda, \mu, \gamma) \in B_{1}$ and $\psi((\lambda, \mu, \gamma)) = (\tilde{\lambda}, \tilde{\mu}, \tilde{\gamma}) \in B_{2}$, $\ell(\mu)-\ell(\lambda)-1 = (\ell(\tilde{\mu})-\ell(\tilde{\lambda})-1)+1$.  Note further that, if $(\lambda, \mu, \gamma) \notin B_{1} \sqcup B_{2}$ and $\psi((\lambda, \mu, \gamma)) = (\tilde{\lambda}, \tilde{\mu}, \tilde{\gamma})$,  $\ell(\mu)-\ell(\lambda)-1 = \ell(\tilde{\mu}) - \ell(\tilde{\lambda})-1$. Since
\begin{eqnarray*} 
&& \sum_{(\lambda, \mu, \gamma) \in S} (-1)^{\ell(\lambda)}
 a^{\ell(\mu)-\ell(\lambda)-1}b^{\ell(\lambda)+\ell(\gamma)}q^{|\lambda|+|\mu|+|\gamma|} = a^{-1}\sum_{n=0}^{\infty} (-1)^{n} q^{\binom{n+1}{2}}b^{n}a^{-n} \prod_{k=1}^{\infty} \frac{1}{(1-aq^{k})(1-bq^{k})} \\
 && \sum_{(\lambda, \mu, \gamma) \in B_{1}} (-1)^{\ell(\lambda)}
 a^{\ell(\mu)-\ell(\lambda)-1}b^{\ell(\lambda)+\ell(\gamma)}q^{|\lambda|+|\mu|+|\gamma|} = a^{-1}\left(\sum_{n=1}^{\infty} \frac{(-1)^{n} q^{\binom{n+1}{2}} b^{n}a^{-n}} {(1-bq)(1-bq^{2}) \cdots (1-bq^{n})} \right)  \\
 && \sum_{(\lambda, \mu, \gamma) \in Fix(\psi)} (-1)^{\ell(\lambda)}
 a^{\ell(\mu)-\ell(\lambda)}b^{\ell(\lambda)+\ell(\gamma)}q^{|\lambda|+|\mu|+|\gamma|}  = a^{-1} + \sum_{n=1}^{\infty} \frac{q^{n}}{(1-aq) \cdots (1-aq^{n})(1-bq) \cdots (1-bq^{n})}\end{eqnarray*}
we have, using the above argument 
\begin{eqnarray*} && a^{-1}\sum_{n=0}^{\infty} (-1)^{n} q^{\binom{n+1}{2}}b^{n}a^{-n} \prod_{k=1}^{\infty} \frac{1}{(1-aq^{k})(1-bq^{k})} - (1-a) a^{-1}\left(\sum_{n=1}^{\infty} \frac{(-1)^{n} q^{\binom{n+1}{2}} b^{n}a^{-n}} {(1-bq)(1-bq^{2}) \cdots (1-bq^{n})} \right) \\ && = a^{-1} + \sum_{n=1}^{\infty} \frac{q^{n}}{(1-aq)(1-aq^{2}) \cdots (1-aq^{n})(1-bq)(1-bq^{2}) \cdots (1-bq^{n})}\end{eqnarray*}
or
\begin{eqnarray*} 
1+ \sum_{n=1}^{\infty} \frac{q^{n}}{(1-aq)(1-aq^{2}) \cdots (1-aq^{n})(1-bq)(1-bq^{2}) \cdots (1-bq^{n})}   \\
 = (1-a^{-1}) \left(1+ \sum_{n=1}^{\infty} \frac{(-1)^{n} q^{\binom{n+1}{2}} b^{n}a^{-n}} {(1-bq)(1-bq^{2}) \cdots (1-bq^{n})} \right)  \nonumber \\ 
+ a^{-1}\sum_{n=0}^{\infty} (-1)^{n} q^{\binom{n+1}{2}}b^{n}a^{-n} \prod_{k=1}^{\infty} \frac{1}{(1-aq^{k})(1-bq^{k})}  \nonumber
\end{eqnarray*}
  after a bit of algebraic manipulation of the $q^{0}$ terms (our one, hopefully minor, deviation from the equation as written).  This suffices to prove (\ref{ramorigin}) using an involution-principle-based combinatorial argument. 
\section{Rewriting Ramanujan's Identity}

Having proven (\ref{ramorigin}) using the involution principle, we will now show that a series of simple algebraic manipulations applied to (\ref{ramorigin}) gives an identity that suggests a natural combinatorial interpretation and bijective proof.  We begin with (\ref{ramorigin}) as written:
\begin{eqnarray*} 
1+ \sum_{n=1}^{\infty} \frac{q^{n}}{(1-aq)(1-aq^{2}) \cdots (1-aq^{n})(1-bq)(1-bq^{2}) \cdots (1-bq^{n})}   \\
 = (1-a^{-1}) \left(1+ \sum_{n=1}^{\infty} \frac{(-1)^{n} q^{\binom{n+1}{2}} b^{n}a^{-n}} {(1-bq)(1-bq^{2}) \cdots (1-bq^{n})} \right)  \nonumber \\ 
+ a^{-1}\sum_{n=0}^{\infty} (-1)^{n} q^{\binom{n+1}{2}}b^{n}a^{-n} \prod_{k=1}^{\infty} \frac{1}{(1-aq^{k})(1-bq^{k})}  \nonumber
\end{eqnarray*}
Making the substitution $c = -b / a$ gives:
\begin{eqnarray*} 1+ \sum_{n=1}^{\infty} \frac{q^{n}}{(1-aq)(1-aq^{2}) \cdots (1-aq^{n})(1+acq)(1+acq^{2}) \cdots (1+acq^{n})}  \\
= (1-a^{-1}) \left(1+ \sum_{n=1}^{\infty} \frac{q^{\binom{n+1}{2}} c^{n}} {(1+acq)(1+acq^{2}) \cdots (1+acq^{n})} \right) \\
+ a^{-1}\sum_{n=0}^{\infty} q^{\binom{n+1}{2}}c^{n} \prod_{k=1}^{\infty} \frac{1}{(1-aq^{k})(1+acq^{k})} \end{eqnarray*}
Multiplying both sides of the result by $a\prod_{k=1}^{\infty}(1+acq^{k})$ gives:
\begin{eqnarray*} a\prod_{k=1}^{\infty}(1+acq^{k}) + \sum_{n=1}^{\infty} \frac{aq^{n}\prod_{k=n+1}^{\infty}(1+acq^{k})}{(1-aq)(1-aq^{2}) \cdots (1-aq^{n})} \\
 = (a-1) \left(\prod_{k=1}^{\infty}(1+acq^{k}) + \sum_{n=1}^{\infty} q^{\binom{n+1}{2}} c^{n}\prod_{k=n+1}^{\infty}(1+acq^{k}) \right) \\
 + \sum_{n=0}^{\infty}q^{\binom{n+1}{2}}c^{n} \prod_{k=1}^{\infty} \frac{1}{(1-aq^{k})} \end{eqnarray*}

Rearranging terms and dividing both sides of the result by $1-a$ gives 
\begin{eqnarray*} \frac{\prod_{k=1}^{\infty}(1+acq^{k})}{1-a} &+& \sum_{n=1}^{\infty} \frac{aq^{n}\prod_{k=n+1}^{\infty}(1+acq^{k})}{(1-a)(1-aq)(1-aq^{2}) \cdots (1-aq^{n})} + \sum_{n=1}^{\infty} q^{\binom{n+1}{2}} c^{n}\prod_{k=n+1}^{\infty}(1+acq^{k})  \\
 &=& \left(\sum_{n=0}^{\infty}q^{\binom{n+1}{2}}c^{n}\right) \left(\prod_{k=0}^{\infty} \frac{1}{(1-aq^{k})} \right)  \end{eqnarray*}
 
Finally, we can combine the first two terms of the left-hand-side, giving the rewritten form of (\ref{ramorigin}) that we will proceed to prove bijectively:

\begin{eqnarray} \label{rambij} \sum_{n=0}^{\infty} \frac{aq^{n}\prod_{k=n+1}^{\infty}(1+acq^{k})}{(1-a)(1-aq)(1-aq^{2}) \cdots (1-aq^{n})} + \sum_{n=1}^{\infty} q^{\binom{n+1}{2}} c^{n}\prod_{k=n+1}^{\infty}(1+acq^{k})  
 = \left(\sum_{n=0}^{\infty}q^{\binom{n+1}{2}}c^{n}\right) \left(\prod_{k=0}^{\infty} \frac{1}{(1-aq^{k})} \right)  \end{eqnarray}
\section{Interpreting and Proving The Rewritten Ramanujan's Identity Bijectively}
Let $R$ be the set of ordered pairs of $(\lambda, \mu)$, such that $\lambda$ is a (possibly empty) ``triangular'' partition and $\mu$ is a (possibly empty) partition which may include parts of size zero. Note that:
\begin{eqnarray} \label{ramgenerating} \sum_{(\lambda, \mu) \in R} a^{\ell(\mu)}c^{\ell(\lambda)}q^{|\lambda|+|\mu|} =
 \left(\sum_{n=0}^{\infty}q^{\binom{n+1}{2}}c^{n}\right) \left(\prod_{k=0}^{\infty} \frac{1}{(1-aq^{k})} \right) \end{eqnarray}
 
Let $A_{1}$ be the set of ordered pairs of partitions $(X, Y)$ such that $X$ is a partition with distinct positive parts, $Y$ is a partition which may include parts of size zero, and $a(Y) < s(X)$.  Let $A_{2}$ be the set of ordered pairs of partitions $(X, Y)$ such that $X$ is a partition with distinct positive parts, $Y$ is a (nonempty)``triangular'' partition, i.e., $Y = (\ell(Y), \ell(Y)-1, \ldots, 1)$, and $a(Y) < s(X)$.  Let $A = A_{1} \sqcup A_{2}$.  

Note that
 \begin{eqnarray*}  
 \sum_{(X, Y) \in A_{1}} a^{\ell(X)+\ell(Y)}c^{\ell(X)} q^{|X|+|Y|} &=& \sum_{n=0}^{\infty}  \frac{aq^{n}\prod_{k=n+1}^{\infty}(1+acq^{k})}{(1-a)(1-aq)(1-aq^{2}) \cdots (1-aq^{n})} \\
  \sum_{(X, Y) \in A_{2}} a^{\ell(X)}c^{\ell(X)+\ell(Y)} q^{|X|+|Y|} &=& \sum_{n=1}^{\infty} q^{\binom{n+1}{2}} c^{n}\prod_{k=n+1}^{\infty}(1+acq^{k}) \end{eqnarray*} 
 
 To prove (\ref{rambij}) bijectively, it therefore suffices to define a bijection $\phi$ between $R$ and $A$, such that, if $\phi((\lambda, \mu)) = (X, Y)$:
 \begin{enumerate}
 \item $|\lambda| + |\mu| = |X|+|Y|$
 \item $\ell(X) + \ell(Y) = \ell(\mu)$ and $\ell(X) = \ell(\lambda)$, if $(X, Y) \in A_{1}$
 \item $\ell(X) + \ell(Y)  = \ell(\lambda)$ and $\ell(X) = \ell(\mu)$, if $(X, Y) \in A_{2}$.  
 \end{enumerate}
 
 In general, for any partitions $A = (A_{1}, A_{2}, \ldots)$ and $B = (B_{1}, B_{2}, \ldots)$, let $A+B = (A_{1}+B_{1}, A_{2}+B_{2}, \ldots)$, where $\ell(A+B) = \max(\ell(A), \ell(B))$.  For example, $(3, 2, 1, 1) + (6, 6, 5) = (9, 8, 6, 1)$.  Note that $|A+B| = |A|+|B|$.  
 
 Define the map $\phi$ on $R$ as follows: For $(\lambda, \mu) \in R$, if $k = \min(\ell(\lambda), \ell(\mu))$, $p = \max(\ell(\lambda), \ell(\mu))$ and $\lambda+\mu = ((\lambda+\mu)_{1}, (\lambda+\mu)_{2}, \ldots, (\lambda+\mu)_{p})$, let $\phi((\lambda, \mu)) = (X, Y)$, where $X = ((\lambda+\mu)_{1}, (\lambda+\mu)_{2}, \ldots, (\lambda+\mu)_{k})$ and $Y = ((\lambda+\mu)_{k+1}, (\lambda+\mu)_{k+2}, \ldots, (\lambda+\mu)_{p})$.  

Note that, given $(\lambda, \mu) \in R$, there are two possibilities:
\begin{enumerate}
\item $\ell(\lambda) \leq \ell(\mu)$.  In this case, $\phi((\lambda, \mu)) = (X, Y)$, where $X = (\mu_{1}+\ell(\lambda), \mu_{2}+\ell(\lambda)-1, \ldots, \mu_{\ell(\lambda)}+1)$ and $Y = (\mu_{\ell(\lambda)+1}, \mu_{\ell(\lambda)+2}, \ldots, \mu_{\ell(\mu)})$.  Then $a(Y) = \mu_{\ell(\lambda)+1} < \mu_{\ell(\lambda)}+1 = s(X)$, so $(X, Y) \in A_{1}$.  Further, $|\lambda| + |\mu| = |X|+|Y|$, $\ell(X) + \ell(Y) = \ell(\mu)$, and $\ell(X) = \ell(\lambda)$.
\item $\ell(\lambda) > \ell(\mu)$.  In this case, $\phi((\lambda, \mu)) = (X, Y)$, where $X = (\mu_{1}+\ell(\lambda), \mu_{2}+\ell(\lambda)-1, \ldots, \mu_{\ell(\mu)}+\ell(\lambda)-\ell(\mu)+1)$ and $Y =(\ell(\lambda)-\ell(\mu), \ell(\lambda)-\ell(\mu)-1, \ldots, 1)$.  Then $a(Y) = \ell(\lambda)-\ell(\mu) < \mu_{\ell(\mu)}+\ell(\lambda)-\ell(\mu)+1 = s(X)$, so $(X, Y) \in A_{2}$.  Further, $|\lambda|+|\mu| = |X|+|Y|$, $\ell(X) +\ell(Y) = \ell(\lambda)$, and $\ell(X) = \ell(\mu)$.
\end{enumerate} 
It remains to prove $\phi$ is a bijection between $R$ and $A$.  Define $\phi^{-1}: A \rightarrow R$ as follows.  Given $(X, Y) \in A$,  with $X = (X_{1}, X_{2}, \ldots, X_{\ell(X)})$ and $Y = (Y_{1}, Y_{2}, \ldots, Y_{\ell(Y)})$:
\begin{enumerate}
\item If $(X, Y) \in A_{1}$, let $\phi^{-1}((X, Y)) = (\lambda, \mu)$, where $\lambda = (\ell(X), \ell(X)-1, \ldots, 1)$ and $\mu = (X_{1}-\ell(X), X_{2}-\ell(X)+1, \ldots, X_{\ell(X)}-1, Y_{1}, Y_{2}, \ldots, Y_{\ell(Y)})$.  Note that $\lambda$ is a ``triangular'' partition, $\mu$ is a partition which may include parts of size zero, and $\ell(\mu) \geq \ell(\lambda)$.  Note also that $|\lambda|+|\mu| = |X|+|Y|$, $\ell(X) + \ell(Y) = \ell(\mu)$, and $\ell(X) = \ell(\lambda)$.     
  
\item If $(X, Y) \in A_{2}$.  In this case, let $\phi^{-1}((X, Y)) = (\lambda, \mu)$, where $\lambda = (\ell(X)+\ell(Y), \ell(X)+\ell(Y)-1, \ldots, 1)$ and $\mu = (X_{1}-\ell(X)-\ell(Y), X_{2} -\ell(X)-\ell(Y)+1, \ldots, X_{\ell(X)} - \ell(Y)-1)$.  Then $\lambda$ is a ``triangular'' partition, $\mu$ is a partition which may include parts of size zero, and $\ell(\lambda) > \ell(\mu)$.  Note also that $|\lambda|+|\mu| = |X|+|Y|$, $\ell(\lambda) = \ell(X) + \ell(Y)$, and $\ell(\mu) = \ell(X)$.  
\end{enumerate}
It should be clear from the parallel division into cases that $\phi$ and $\phi^{-1}$ are in fact inverses of each other, which suffices to prove (\ref{rambij}).  
\subsection{Examples}
Figures $1$ and $2$ give two examples of elements $(\lambda, \mu) \in R$, with $\lambda = (6, 5, 4, 3, 2, 1)$ in each case.  In Figure $1$, $(\lambda, \mu) = ((6, 5, 4, 3, 2, 1), (8, 6, 6, 6, 4, 4, 4, 3, 3, 0, 0, 0, 0,0))$, with 
\[ \phi(\lambda, \mu)) = ((14, 11, 10, 9, 6, 5), (4, 3, 3, 0, 0, 0, 0, 0)) \in A_{1}\].  In Figure $2$, $(\lambda, \mu) = ((6, 5, 4, 3, 2, 1), (8, 8, 0))$, with $\phi((\lambda, \mu)) =  ((14, 13, 4), (3, 2, 1)) \in A_{2}$.  

\begin{figure}[ht]
\includegraphics[width = 10 cm]{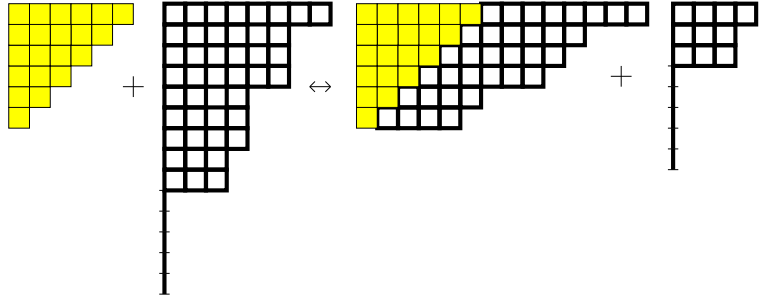}
\caption{$\phi(((6, 5, 4, 3, 2, 1), (8, 6, 6, 6, 4, 4, 4, 3, 3, 0, 0, 0, 0,0))) = ((14, 11, 10, 9, 6, 5), (4, 3, 3, 0, 0, 0, 0, 0))$}
\end{figure}

\begin{figure}[ht]
\includegraphics[width = 10 cm]{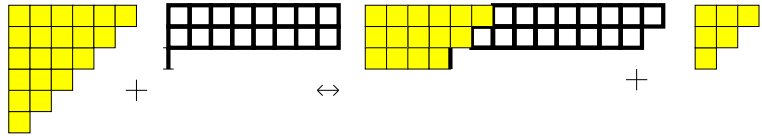}
\caption{$\phi(((6, 5, 4, 3, 2, 1), (8, 8, 0))) = ((14, 13, 4), (3, 2, 1))$}
\end{figure}
\section{Note on the Bijection Used To Prove the Rewritten Identity}
An elementary result in partition theory is that $(\mu_{1}, \mu_{2}, \ldots, \mu_{k}) \leftrightarrow (\mu_{1}+k, \mu_{2}+k-1, \ldots, \mu_{k}+1)$ gives a bijection between partitions of $n$ into $k$ non-negative parts and partitions of $n+\binom{k+1}{2}$ into $k$ distinct positive parts.  Using the above notation, this bijection can be written as $\mu \leftrightarrow \mu + \lambda$, where $\ell(\mu) = \ell(\lambda) = k$ and $\lambda$ is a triangular partition--in other words, the bijection operates by adding, to a partition, a triangular partition of the same length.  Ramanujan's identity (\ref{ramorigin}), when rewritten as (\ref{rambij}), can be interpreted, as we did above, as exploring the result of adding, to a partition, a triangular partition \textit{not necessarily of the same length}.  
\bibliography{mybib}{}

\begin{thebibliography}{1}

\bibitem{Andrewsint}
George~E. Andrews.
\newblock An introduction to {R}amanujan's ``lost'' notebook.
\newblock {\em Amer. Math. Monthly}, 86(2):89--108, 1979.

\bibitem{Andrewspartialtheta}
George~E. Andrews and S.~Ole Warnaar.
\newblock The product of partial theta functions.
\newblock {\em Adv. in Appl. Math.}, 39(1):116--120, 2007.

\bibitem{kimpartialtheta}
Byungchan Kim.
\newblock Combinatorial proofs of certain identities involving partial theta
  functions.
\newblock {\em Inter. J. Number Theory, in press}.

\bibitem{pakinv}
Igor Pak.
\newblock The nature of partition bijections. {I}. {I}nvolutions.
\newblock {\em Adv. in Appl. Math.}, 33(2):263--289, 2004.

\bibitem{paksurvey}
Igor Pak.
\newblock Partition bijections, a survey.
\newblock {\em Ramanujan J.}, 12(1):5--75, 2006.

\end{thebibliography}
\bibliographystyle{plain}
 \end{document}